\documentclass[12pt, a4paper]{article}
\usepackage[T1]{fontenc}
\usepackage[french]{babel}
\usepackage{amsmath,amsthm}
\usepackage{amssymb}
\usepackage{amscd}
\usepackage[latin1]{inputenc}
\usepackage{amsfonts}
\usepackage[matrix,arrow,curve]{xy}
\usepackage{graphicx}

\theoremstyle{plain}
\newtheorem{theo}{Th\'eor\`eme}[section]
\newtheorem{lem}[theo]{Lemme}
\newtheorem{prop}[theo]{Proposition}
\newtheorem{cor}[theo]{Corollaire}

\theoremstyle{definition}

\newtheorem{theosd}[theo]{Th\'eor\`eme}
\newtheorem{propsd}[theo]{Proposition}
\newtheorem{rem}[theo]{Remarque}

\newtheorem{corro}[theo]{Corollaire}

\newcommand{\al}{\ensuremath{\alpha}}

\newcommand{\hh}{\ensuremath{\mathcal{H}}}
\newcommand{\QZl}{\ensuremath{\mathbb Q_l/\mathbb Z_l}}
\newcommand{\QZ}{\ensuremath{\mathbb Q/\mathbb Z}}
\newcommand{\FF}{\ensuremath{\mathbb{F}}}
\newcommand{\ZZ}{\ensuremath{\mathbb Z}}
\newcommand{\ZZl}{\ensuremath{\mathbb Z}_l}
\newcommand{\HH}{\ensuremath{\mathbb H}_{\acute{e}t}}
\newcommand{\nH}{\ensuremath{H_{\mathrm{nr}}}}

\title{Invariants birationnels dans la suite spectrale de Bloch-Ogus}
\author{Alena Pirutka}

\begin{document}

\maketitle

\begin{abstract}
On établit l'invariance birationnelle des groupes $H^{i}(X, \hh^{n+d}(\mu_{l^r}^{\otimes j}))$ pour $X$ une variété projective et lisse, géométriquement intègre, de dimension $n$, définie sur un corps $k$ de dimension cohomologique au plus $d$, $(l, \mathrm{car}.k)=1$.
On obtient aussi un résultat analogue sur un corps fini pour les groupes $H^{i}(X, \hh^{n}(\QZl(j)))$ et on relie un de ces invariants avec le conoyau de l'application classe de cycle
$CH^{n-1}(X)\otimes \ZZl \to H^{2n-2}_{\acute{e}t}(X,\ZZl(n-1))$, ce qui donne une version sur un corps fini d'un résultat de Colliot-Thélène et Voisin \cite{CTV} 3.11 sur le corps des complexes. 
\end{abstract}

\section{Introduction}

Soit $k$ un corps. 
Pour $n$ un entier  inversible sur $k$,  on note $\mu_{n}$ le $k$-sch\'ema en groupes (\'etale) des racines $n$-i\`emes de l'unit\'e. Pour $j$ un entier positif on note $\mu_{n}^{\otimes j}=\mu_{n}\otimes\ldots\otimes\mu_{n}$ ($j$ fois). On pose $\mu_{n}^{\otimes j}=Hom_{k-gr}(\mu_{n}^{\otimes (-j)}, \mathbb Z/n)$ si $j$ est n\'egatif et $\mu_{n}^{\otimes 0}=\mathbb Z/n$.  Ces $k$-sch\'emas en groupes donnent des faisceaux \'etales, not\'es encore $\mu_{n}^{\otimes j}$, sur toute  vari\'et\'e définie sur $k$.  %On note $H^i(X,\mu_n^{\otimes j})$ les groupes de cohomologie \'etale de $X$ \`a valeurs dans $\mu_n^{\otimes j}$. 

Soit $X$ une $k$-variété intègre, projective et lisse. 
On note $\mathcal H^i_X(\mu_n^{\otimes j})$ le faisceau de Zariski sur $X$ associ\'e au pr\'efaisceau $U\mapsto H^i_{\acute{e}t}(U,\mu_n^{\otimes j})$. La conjecture de Gersten, établie par Bloch et Ogus (\cite{BO}) permet de calculer les groupes de cohomologie de ces faisceaux comme les groupes de cohomologie du complexe
\small
\begin{equation*}
0\to H^{i}(k(X),\mu_n^{\otimes j}) \to  \bigoplus_{x\in X^{(1)}}H^{i-1}(\kappa(x),\mu_n^{\otimes( j-1)})\to\ldots\to  \bigoplus_{x\in X^{(r)}}H^{i-r}(\kappa(x),\mu_n^{\otimes (j-r)})\to\ldots
\end{equation*}
\normalsize
où $X^{(r)}$ désigne l'ensemble des points de codimension $r$ de $X$; $\kappa(x)$  est le corps résiduel du point $x$. Les flèches de ce complexe sont induites par des résidus et  le terme $\bigoplus_{x\in X^{(r)}}$ est en degré $r$. 

  On a une suite spectrale  (cf. \cite{BO}) 
\begin{equation}\label{sspbo}
 E_2^{pq}= H^p(X, \mathcal H^q(\mu_n^{\otimes j}))\Rightarrow H^{p+q}_{\acute{e}t}(X,\mu_n^{\otimes j}).
\end{equation}

Les termes $E_2^{0q}$ de cette suite spectrale s'identifient à des groupes de cohomologie non ramifiée $\nH^i(X, \mu_n^{\otimes j})$ qui sont des invariants birationnels des $k$-variétés intègres projectives et lisses (cf. \cite{CT92}). Dans ce texte on s'intéresse à d'autres invariants birationnels dans (\ref{sspbo}).

Pour $k$ un corps de dimension cohomologique au plus $d$ on établie l'invariance birationnelle des groupes $H^{i}(X, \hh^{n+d}(\mu_{l^r}^{\otimes j}))$. Ceci est  fait  dans la section \ref{sib} par un argument d'action des correspondances.  Cette action est fournie par la théorie de cycles de Rost, dont on fait des rappels dans la section \ref{rr}. D'après la conjecture de Kato (\cite{Ka}, \cite{KeS}) les groupes $H^{i}(X, \hh^{n+1}(\mu_{l^r}^{\otimes n}))$, $i<n$, sont nuls pour   $X$ une variété projective et lisse, géométriquement intègre, de dimension $n$, définie sur un corps fini $\FF$, avec $r$ premier à la caractéristique de $\FF$. Dans ce cas, on établit l'invariance birationnelle des groupes $H^{i}(X, \hh^{n}(\QZl(j)))$, cf.  \ref{ibf}. Dans la section \ref{luc},  on montre  que le quotient du groupe $H^{n-3}(X,\mathcal H^n(\QZl(n-1)))$ par son sous-groupe divisible maximal est isomorphe au  groupe de torsion du  conoyau  de l'application classe de cycle
$CH^{n-1}(X)\otimes \ZZl \to H^{2n-2}_{\acute{e}t}(X,\ZZl(n-1))$.  %Pour ce faire, on a besoin de la conjecture de Kato pour des twists $j$ autres que $n-1$. En utilisant l'argument de  Kahn \cite{Ka93}, on l'établie dans la section \ref{rck}.% (cf. aussi appendice).

\paragraph*{Remerciements.} Je voudrais remercier mon directeur de thèse, Jean-Louis Colliot-Thélène, de m'avoir introduit aux problèmes considérés dans ce texte, pour de nombreuses discussions et pour son aide constante.

\section{Rappels sur les modules de cycles de Rost}\label{rr}

 Tous les énoncés de ce paragraphe se trouvent dans les articles de Rost \cite{Ro} et de Déglise \cite{D}.

 Soit $k$ un corps, soit $X$ un $k$-schéma équidimensionnel et soit $\mathcal F(X)$ une classe de corps sur $X$ (corps qui contiennent un corps r\'esiduel d'un point de $X$).% (par exemple, la classe de corps résiduels de points de $X$).

\begin{enumerate}
\item On peut voir un \textbf{module de cycles} comme un foncteur $M:\mathcal F(X)\to Ab$, $M=\coprod M_n$, qui satisfait certains axiomes (existence des analogues de restriction, corestriction, résidus, multiplication par $K_1$ et compatibilités entre ces applications).\\
    \textit{Exemples.} \begin{itemize} \item[$\bullet$] $M_H(F)=\coprod_nH^n(F, D\otimes\mu_r^{\otimes n})$ où $D$ est un $G$-module fini continu, d'exposant $r$. % En particulier, on récupère ici les  $H^n(F, \mu_r^{\otimes j})$;
 \item[$\bullet$] $M_K(F)=\coprod_n K^M_n(F)$. \\Le groupe $K_n^M(F)$ est le $n$-ième  groupe de $K$-théorie de Milnor de $F$. C'est le quotient de $\underbrace {F^*\otimes\ldots\otimes F^*}_{\mbox{n fois}}$ par le sous-groupe engendré par des éléments $a_1\otimes\ldots\otimes a_n$ avec $a_i+a_j=1$ pour $1\leq i<j\leq n$. En particulier, $K^M_0(F)=\mathbb Z$ et $K^M_1(F)=F^*$.
\end{itemize}
\item Un \textbf{accouplement} $M\times M'\to M''$ de modules de cycles est la donnée pour tout $F\in \mathcal F(X)$ de  $M(F)\times M'(F)\to M''(F)$, qui satisfait des propriétés de compatibilité.\\
    \textit{Exemple.} Pour tout module de cycles $M$, la multiplication par $K_1$ donne un accouplement $M_K\times M\to M$.
\item \textbf{Complexes et groupes de Chow}. Les groupes $C^p(X,M):=\coprod\limits_{x\in X^{(p)}} M(\kappa(x))$ forment un complexe $$C(X,M)=[\ldots\to C^p(X,M)\to C^{p+1}(X,M))\to\ldots]$$ et on note $$A^p(X,M)=H^p(C(X,M)).$$
    \textit{Exemples.} \begin{itemize} \item[$\bullet$] D'après la définition, le groupe de Chow $CH^p(X)$ est un facteur direct de $A^p(X,M_K)$.
\item[$\bullet$] Pour $X$ lisse, d'après les r\'esultats de Bloch et Ogus \cite{BO}, on a  $A^p(X,M_H)=\coprod_nH^p(X,\mathcal H^n(D\otimes\mu_r^{\otimes n}))$, où $\mathcal H^n(D\otimes\mu_r^{\otimes n})$ d\'esigne le faisceau de Zariski sur $X$ associ\'e au pr\'efaisceau $U\mapsto H^j(U,D\otimes\mu_r^{\otimes n})$.
\end{itemize}
\item \textbf{Fonctorialités}.\begin{enumerate} \item Pour $f:Y\to X$ propre, on a $f_*: A^p(Y,M)\to A^{p-d}(X,M)$ où $d=dim\,Y-dim\,X$ est la dimension relative de $f$ (pour simplifier, pour $X$ et $Y$ intègres). Ce morphisme est induit  par un morphisme de complexes de bidegré $(d,0)$, correspondant respectivement à la graduation de $C$ et de $M$ (cf. \cite{D} 1.3).  Il est  fonctoriel, c'est-à-dire, pour $g : Z\to Y$  propre, on a $(f\circ g)_*=f_*\circ g_*$ (\cite{Ro} 4.1). \\ \textit{Exemples.}\begin{itemize} \item[$\bullet$] Pour $Y$ de codimension $1$ dans $X$ et $f:Y\to X$ une immersion fermée on a $$H^p(Y,\mathcal H^n(D\otimes\mu_r^{\otimes n}))\to H^{p+1}(X,\mathcal H^{n+1}(D\otimes\mu_r^{\otimes {(n+1)}})).$$
    \item [$\bullet$]Pour $Y$ propre et $f:X\times Y\to X$ une projection, on a $$H^p(X\times Y,\mathcal H^n(D\otimes\mu_r^{\otimes n}))\to H^{p-d}(X,\mathcal H^{n-d}(D\otimes\mu_r^{\otimes {(n-d)}}))$$ où $d=dim\,Y$.
    \end{itemize}
    \item  Pour $f:Y\to X$ un morphisme plat, on dispose d'un morphisme $f^*: A^p(X,M)\to A^{p}(Y,M)$,  induit  par un morphisme de complexes  de bidegré $(0,0)$. Dans le cas où $X$ et $Y$ sont lisses, par la construction plus difficile (\cite{Ro} 12, \cite{D} 3.18) on dispose d'un morphisme $f^*: A^p(X,M)\to A^{p}(Y,M)$ pour $f:Y\to X$ quelconque. Pour $g:Z\to Y$ avec $Z$ lisse, on a $(f\circ g)^*=g^*\circ f^*$ (\cite{Ro} 12.1).
        %\\ \textit{Exemple.} Cela donne $f^*: H^p(X,\mathcal H^n(D\otimes\mu_r^{\otimes n}))\to H^{p}(Y,\mathcal H^{n}(D\otimes\mu_r^{\otimes {n}}))$.
\end{enumerate}
\item \textbf{Localisation}. Pour $Y\stackrel{i}{\to}X$ un fermé purement de codimension $d$ et $U=X\setminus Y\stackrel{j}{\hookrightarrow} X$ son complémentaire, on a une longue suite exacte de localisation\footnote{sans hypothèses de lissité} (\cite{Ro} 5, p.356)$$\ldots\stackrel{\partial}{\to} A^{p-d}(Y,M)\stackrel{i_*}{\to} A^p(X,M)\stackrel{j^*}{\to} A^p(U,M)\stackrel{\partial}{\to}A^{p-d+1}(Y,M)\to\ldots$$
\item \textbf{Invariance homotopique}. Pour $\pi:\mathbb A^n_X\to X$ la projection naturelle, l'application $\pi^*:A^p(X,M)\to A^{p}(\mathbb A^n_X, M)$ est un isomorphisme\footnote{sans hypothèses de lissité} (\cite{Ro} 8.6).
\item \textbf{Cup-produits}. Pour un accouplement $N\times M\to M$ de modules de cycles, on a des produits $$\times: C^p(Y,N)\times C^q(Z,M)\to C^{p+q}(Y\times Z, M).$$ Pour $X$ lisse cela donne $$\cup:  A^*(X,N)\times A^*(X,M)\to A^*(X\times X, M)\stackrel{\Delta_X^*}{\to} A^*(X,M)$$ où $\Delta_X$ est la diagonale.
    On a ici une formule de projection (\cite{D} 5.9(3)) pour $X, Y$ lisses et $f:Y\to X$ propre $$f_*(x\cup f^*y)=f_*x\cup y.$$
    \textit{Exemple.} En utilisant que $CH^p(X)$ est un facteur direct de $A^p(X,M_K)$, pour $X$ propre et lisse, l'accouplement $M_K\times M\to M$ donne $$CH^p(X)\times A^q(X,M)\to A^{p+q}(X, M).$$
    Pour $M=M_H$ cela donne
    $$CH^p(X)\times H^q(X,\mathcal H^n(D\otimes\mu_r^{\otimes n}))\to H^{p+q}(X,\mathcal H^{n+p}(D\otimes\mu_r^{\otimes(n+p)})).$$
\item \textbf{Action de correspondances}. (cf. \cite{D} $\S 6$) \'Etant donné une correspondance $\al\in CH^p(X\times Y)$ avec $X$ et $Y$ propres et lisses, on définit une application $\alpha_*:A^q(X,M)\to A^{p+q-dimX}(Y,M)$ par la formule usuelle
    $$ \alpha_*(x)=p_{Y, *}(\alpha\cup p_X^*x),$$ où $p_X$ (resp. $p_Y$) est une projection de $X\times Y$ sur le premier (resp. sur le deuxième) facteur. Au moins pour les correspondances finies, i.e. pour les correspondaces $\al\in CH^p(X\times Y)$  qui sont finies et surjectives sur une composante connexe de $X$, cette action est compatible avec la composition de correspondances, d'après  \cite{D} 6.5. Ce cas de correspondances finies nous suffira pour la suite. Dans le cas général, on peut  voir qu'il suffit de vérifier  les trois propriétés de \cite{CTV} 9.3; on les trouve dans \cite{D} 5.9, puis dans \cite{Ro} 12.5.
\end{enumerate}

\section{Application aux invariants birationnels}\label{sib}

Soit $k$ un corps et soit $X$ une $k$-variété projective et lisse, géométriquement intègre. Soit $l$ un nombre premier, $(l, \mathrm{car}\,.k)=1$. On peut voir les groupes $H^i(X, \hh^{m}(\mu_{l^r}^{\otimes j}))$ comme les groupes de Chow, associés à des modules de cycles de Rost. On dispose  ainsi d'une action de correspondances sur ces groupes.

\lem\label{lnul}{Soit $k$ un corps   et soit $X$ une $k$-variété projective et lisse, géométriquement intègre, de dimension $n$. Soit $Z\subset X\times X$ un cycle premier de dimension $n$, dont l'image par la première projection est contenue dans une sous-variété fermée propre $D\subset X$. Soient $i,m\geq 0$ et $j$ des entiers. Soit $l$ un nombre premier, $(l, \mathrm{car}\,. k)=1$. Supposons que le groupe $H^i(V, \hh^{m}(\mu_{l^r}^{\otimes j}))$ est nul pour toute $k$-variété projective et lisse $V$, de dimension $n-1$. Alors pour toute classe $\alpha\in H^i(X, \hh^{m}(\mu_{l^r}^{\otimes j}))$ on a $[Z]_*(\alpha)=0$.}
\proof{On peut supposer que $D$ est un diviseur sur $X$.

%Supposons d'abord que $k$ est  parfait.
 D'après le théorème de de Jong, amélioré par Gabber (cf.\cite{I}), il existe un morphisme propre, génériquement fini  $f:X'\to X$ dont le degré $d$ est premier à $l$, tel que $X'$ est lisse et $f^{-1}(D)_{\mathrm{red}}$ est à croisements normaux, en particulier, ses composantes irréductibles sont lisses.

On a le diagramme commutatif suivant, où l'on note $X_2=X\times X, X_2'=X'\times X'$,  $\hh^{m}=\hh^{m}(\mu_{l^r}^{\otimes j})$ et $\hh^{m+n}=\hh^{m+n}(\mu_{l^r}^{\otimes (n+j)})$:
\footnotesize
$$\xymatrix{ CH^n(X_2')\times H^i(X', \hh^{m}) \ar[r]^{id\times pr_1^*}\ar@<-7ex>[d]^{f_*}& CH^n(X_2')\times H^i(X_2', \hh^{m})\ar[r]^(.55){\cup}\ar@<-7ex>[d]^{f_*} & H^{n+i}(X_2', \hh^{m+n})\ar[r]^{ pr_{2,*}}\ar[d]^{f_*}&H^{i}(X', \hh^{m})\ar[d]^{f_*}\\
 CH^n(X_2)\times H^i(X, \hh^{m}) \ar[r]^{id\times pr_1^*} \ar@<-6ex>[u]^{f^*}&CH^n(X_2)\times H^i(X_2, \hh^{m})\ar[r]^(.55){\cup} \ar@<-6ex>[u]^{f^*}& H^{n+i}(X_2, \hh^{m+n})\ar[r]^{ pr_{2,*}}&H^{i}(X, \hh^{m})}$$
\normalsize

Le carré au milieu de ce diagramme commute d'après la formule de projection. Soit $\al\in H^i(X, \hh^{m})$. Puisque $f_*f^*[Z]=d[Z]$ dans $ CH^n(X_2)$, on a
\begin{multline*}d[Z]_*(\al)=pr_{2,X,*}(d[Z]\cup pr_{1,X}^*\al)=pr_{2,X,*}(f_*(f^*[Z]\cup f^*pr_{1,X}^*\al))=\\=f_*(pr_{2,X',*}(f^*[Z]\cup pr_{1,X'}^*f^*\al))=f_*([f^*Z]_*(f^*\al)).
\end{multline*}
 Puisque $d$ est premier à $l$, il suffit donc de montrer que $[f^*Z]_*( f^*\al)=0$. On peut le faire pour chaque composante irréductible séparément. On peut donc supposer que  $Z':=f^*Z$  est inclus dans  $D'\times X'$  avec $D'$ lisse. On écrit  $\imath$ pour les inclusions $D'\hookrightarrow X'$ et  $D'\times X' \hookrightarrow X\times X'$.

On a le diagramme

\footnotesize
$$\xymatrix{ CH^{n-1}(D'\times X')\times H^i(D', \hh^{m})\ar[r]^{id\times pr_1^*}\ar@<-7ex>[d]^{\imath_*} & CH^{n-1}(D'\times X')\times H^i(D'\times X', \hh^{m})\ar[r]^(.7){pr_{2,*}\circ\cup}\ar@<-7ex>[d]^{\imath_*} & H^{i}(X', \hh^{m})\ar[d]\\
 CH^n(X'\times X')\times H^i(X', \hh^{m})\ar[r]^{id\times pr_1^*} \ar@<-6ex>[u]^{\imath^*}& CH^n(X'\times X')\times H^i(X'\times X', \hh^{m})\ar[r]^(.7){pr_{2,*}\circ\cup} \ar@<-6ex>[u]^{\imath^*}& H^{i}(X', \hh^{m}).}$$
\normalsize

On a donc que l'action de $[Z']$ se factorise par $H^{i}(D', \hh^{m})$, groupe qui est nul d'après  l'hypothèse. Cela termine la preuve du lemme.\qed% pour $k$ parfait. Puisque $(l,\mathrm{car.}k)=1$, dans le cas général, on passe à une clôture parfaite de $k$ et on déduit le résultat par un argument de corestriction.
%Dans le cas général, soit $K'$  une clôture parfaite de $K$. On a le diagramme commutatif

%\begin{center}
%\xymatrix{
% H^i(X, \hh^{m}(\mu_{l^r}^{\otimes j}))  \ar[r]^{[Z]_*}\ar[d] & H^i(X, \hh^{m}(\mu_{l^r}^{\otimes j}))\ar[d]&\\
% H^i(X_K', \hh^{m}_{X_{K'}}(\mu_{l^r}^{\otimes j})) \ar[r]^{[Z]_*} &H^i(X_K', \hh^{m}_{X_{K'}}(\mu_{l^r}^{\otimes j})).&
%}
%\end{center}
% Ainsi, il existe une extension finie $K''/K$ dont le degré est une puissance de ${\mathrm{car}.K}$, telle que $Res_{K''/K}(\beta)$ provient d'un élément de $ H^3(K'', \QZl(2))$. Puisque $l\neq{\mathrm{car}.K}$, on obtient le résultat par un argument de  corestriction.
%\qed\\}

\rem{\begin{itemize}
\item[(i)] Ce lemme s'applique aussi à des groupes $H^i(X, \hh^{m}(\QZl(j)))$.
\item[(ii)] Plus généralement, le lemme s'applique à des groupes $A^i(X,M_m)$ pour un module de cycles $M=\coprod M_m$, sous les hypothèses que $A^i(X,M_m)$ est de torsion $l$-primaire avec $(l,\mathrm{car}.k)=1$ et que  $A^i(V,M_m)=0$ pour toute $k$-variété $V$ projective et lisse de dimension $n-1$.\\
\end{itemize}}

\theo\label{inv1}{Pour $k$ un corps de dimension cohomologique $cd\,k\leq d$, les groupes $H^i(X, \hh^{n+d}(\mu_{l^r}^{\otimes j}))$, où $ j\in \mathbb Z$ et $(l,\mathrm{car.}k)=1$, sont des invariants birationnels des $k$-variétés projectives et lisses, géométriquement intègres,  de dimension $n$. De plus,
\begin{itemize}
\item [(i)] pour $X/k$ lisse, $H^{i}(X, \hh^{n+d}(\mu_{l^r}^{\otimes j}))\simeq H^{i+N}(\mathbb P^N_X, \hh^{n+N+d}(\mu_{l^r}^{\otimes (j+N)}))$;
\item [(ii)] on a $H^i(\mathbb P^N_k, \hh^{N+d}(\mu_{l^r}^{\otimes j}))=0$ pour $i<N$ et $H^N(\mathbb P^N_k, \hh^{N+d}(\mu_{l^r}^{\otimes j}))\simeq H^d(k, \mu_{l^r}^{\otimes (j-N)})$;
\item [(iii)] si $X$ et $Y$ sont des $k$-variétés projectives et lisses, géométriquement intègres, de dimensions respectives $n_X$ et $n_Y$, stablement birationnelles, i.e. telles que $X\times \mathbb P^{N_X}_k$ et $Y\times \mathbb P^{N_Y}_k$ sont birationnelles, alors $$H^i(X, \hh^{n_X+d}(\mu_{l^r}^{\otimes j}))\simeq H^{i+n_Y-n_X}(Y, \hh^{n_Y+d}(\mu_{l^r}^{\otimes (j+n_Y-n_X)})), \quad i\geq 0.$$
\end{itemize}}
\proof{On raisonne  comme dans \cite{CTV} 3.4. Soient $X,Y$ deux $k$-variétés géométriquement int\`egres, projectives et lisses, et soit $\phi : X\dashrightarrow Y$ une application birationnelle. On note $n=\mathrm{dim}X=\mathrm{dim}Y$.  L'adhérence $\Gamma_{\phi}$ du graphe de $\phi$ définit une correspondance $[\Gamma_{\phi}]\in  CH^n(X\times Y)$. On note $\Gamma_{\phi^{-1}}\subset Y\times X$ le graphe de $\phi^{-1}$. D'après \cite{CTV} 3.5, le composé $[\Gamma_{\phi^{-1}}]\circ[\Gamma_{\phi}]\in CH^n(X\times X)$ se décompose comme $[\Gamma_{\phi^{-1}}]\circ[\Gamma_{\phi}]=\Delta_X+W$
où $\Delta_X$ est la diagonale de $X$ et $W$ est un cycle sur $X\times X$ supporté sur $D\times X$, $D$ étant une sous-variété fermée propre de $X$.

 Les correspondances $[\Gamma_{\phi}]$ et $[\Gamma_{\phi^{-1}}]$ permettent de définir des applications entre les groupes $H^i(X, \hh^{n+d}(\mu_{l^r}^{\otimes j}))$ et $H^i(Y, \hh^{n+d}(\mu_{l^r}^{\otimes j}))$ (cf. section \ref{rr}). D'après Bloch-Ogus et l'hypothèse $cd\,k\leq d$,  $H^i(V, \hh^{n+d}(\mu_{l^r}^{\otimes j}))=0$ pour toute $k$-variété projective et lisse $V$, de dimension $n-1$. Le lemme ci-dessus, appliqué à chaque composante irréductible de  $W$ (resp. à chaque composante de $[\Gamma_{\phi}]\circ[\Gamma_{\phi^{-1}}]-\Delta_Y$), montre alors que les applications composées induites sont des  identités, d'où le premier énoncé du théorème.

 D'après ce qui précède,  pour établir $(i)$, il suffit de montrer que $H^i(X, \hh^{n+d}(\mu_{l^r}^{\otimes j}))\simeq H^{i+1}(X\times\mathbb P^1, \hh^{n+d+1}(\mu_{l^r}^{\otimes (j+1)}))$. On écrit la suite de localisation pour $Y=X\times \{\infty\}\subset X\times\mathbb P^1$ : \begin{multline*}H^{i}(X\times\mathbb A^1, \hh^{n+d+1}(\mu_{l^r}^{\otimes (j+1)})){\to} H^i(X, \hh^{n+d}(\mu_{l^r}^{\otimes j})) {\to} H^{i+1}(X\times\mathbb P^1, \hh^{n+d+1}(\mu_{l^r}^{\otimes (j+1)}))\to\\\to H^{i+1}(X\times\mathbb A^1, \hh^{n+d+1}(\mu_{l^r}^{\otimes (j+1)})).\end{multline*}
\normalsize

 On a $$H^{s}(X\times\mathbb A^1, \hh^{n+d+1}(\mu_{l^r}^{\otimes (j+1)}))\stackrel{pr_X^*}{\simeq} H^{s}(X, \hh^{n+d+1}(\mu_{l^r}^{\otimes (j+1)}))=0, \;s\geq 0$$  d'après Bloch-Ogus ($dim\,X=n$ et $cd\,k\leq d$). On applique cela à $s=i \text{ et }i+1$ et on obtient l'énoncé  $(i)$.

Les énoncés $(ii)$ et $(iii)$ résultent de $(i)$, ce qui termine la preuve du théorème. \qed\\}

\section{Invariants sur un corps fini}

%Soit  $\FF$ un corps fini de caractéristique $p$.

\subsection{Rappels sur la conjecture de Kato}\label{rck}

Soit $X$ un schéma de type fini sur $\FF$. Soit $n$ un entier, $(n,\mathrm{car}.\FF)=1$. Dans \cite{Ka}, Kato a introduit  le  complexe $KC(X,\mathbb Z/n \mathbb Z)$ suivant  :
\begin{multline*}
\ldots\to \bigoplus_{x\in X_{(i)}}H^{i+1}(\kappa(x),\mu_n^{\otimes i})\to  \bigoplus_{x\in X_{(i-1)}}H^{i}(\kappa(x),\mu_n^{\otimes (i-1)})\to\\
\ldots\to \bigoplus_{x\in X_{(1)}}H^2(\kappa(x),\mu_n))\to \bigoplus_{x\in X_{(0)}}H^1(\kappa(x),\mathbb Z/n \mathbb Z))
\end{multline*}
où le terme $\bigoplus_{x\in X_{(i)}}$ est en degré $i$. On note $KH_i(X,\mathbb Z/n \mathbb Z)$ le $i$-ème groupe d'homologie de ce complexe.

Pour $l$ un nombre premier différent de la caractéristique de $\FF$, on pose  $KC(X,\QZl)=\varinjlim KC(X,\mathbb Z/l^m\mathbb Z)$ et on note $KH_i(X,\QZl)$ son $i$-ème groupe d'homologie.

Pour $X$ projective et lisse, géométriquement intègre, Kato \cite{Ka} a conjecturé l'exactitude de ces complexes, sauf au terme $0$. Cette conjecture a été récemment démontrée  par Kerz et Saito :
\theosd(Kerz et Saito \cite{KeS}, Thm.0.4)\label{tks} {\textit{Soit $X$ une variété projective et lisse sur $\FF$ et soit $l$ un nombre premier, $(l,\mathrm{car}.\FF)=1$. Alors $KH_i(X,\mathbb Z/l^m \mathbb Z)=0$ pour tout $i>0$ et $m>0$ et
$KH_0(X,\mathbb Z/l^m \mathbb Z)\simeq \mathbb Z/l^m \mathbb Z$.\\}} %{\textit{Soit $X$ une variété projective et lisse sur $\FF$ et soit $l$ un nombre premier, $(l,p)=1$. Alors $KH_i(X,\QZl)=0$ pour tout $i>0$ et $KH_0(X,\QZl)=\QZl$.\\}}

%Sous la conjecture de Bloch-Kato, qui est maintenant établi par Voevodsky \cite{Voe}, on a aussi l'exactitude de $KH_i(X,\mathbb Z/l^m \mathbb Z)$ en degrés positifs :

%\theosd\label{tjs}(\cite{JS}) {\textit{Soit $X$ une variété projective et lisse sur $\FF$ et soit $l$ un nombre premier, $(l,p)=1$. Alors $KH_i(X,\mathbb Z/l^m \mathbb Z)=0$ pour tout $i>0$ et $m>0$ et
%$KH_0(X,\mathbb Z/l^m \mathbb Z)\simeq \mathbb Z/l^m \mathbb Z$.\\}}

D'après la définition, pour $X$ une $\FF$-variété projective et lisse, de dimension $n$, $H^i(X, \hh^{n+1}(\mu_{l^m}^{\otimes n}))=KH_{n-i}(X, \mathbb Z/l^m)$. Ces invariants birationnels sont donc tous nuls, sauf le dernier.

Pour $j\neq n$ et pour les coefficients $\QZl$, un lemme  de  Sato et Saito  permet de voir aussi la nullité des groupes $H^i(X, \hh^{n+1}((\QZl(j)))$, ce qui n'utilise pas la conjecture de Kato. En effet, on a dans ce cas que tous les termes  dans la résolution de Bloch-Ogus du faisceau  $\hh^{n+1}(\QZl(j))$  sont nuls.

\propsd\label{cks}( Kahn \cite{Ka93}; Saito et Sato \cite{SS}, Lemma 2.7) {\textit{ Soit $L$ un corps de degré de transcendance $r$ sur $\FF$. Soit $l$ un nombre premier, $(l,\mathrm{car}.\FF)=1$. Alors $$H^{r+1}(L,\QZl(j))=0$$ pour tout $j\neq r$. En  particulier, pour $X$ une variété projective et lisse sur $\FF$ de dimension $n$ et pour tout $j\neq n$ le faisceau $\hh^{n+1}(\QZl(j))$ est identiquement nul, ainsi que tous les termes de sa résolution de Bloch-Ogus. Ainsi,  $H^i(X, \hh^{n+1}(\QZl(j)))=0$ pour tout $i$ et pour $j\neq n$.\\}

Quant aux coefficients finis, en utilisant la méthode de Bruno Kahn \cite{Ka93}, on établit aussi  la nullité des groupes $H^i(X, \hh^{n+1}(\mu_{l^m}^{\otimes j}))$ pour $j\in \mathbb Z$ et $i<n$ comme une conséquence de la conjecture de Kato.  Cette méthode utilise  également la conjecture de Bloch-Kato. Donnons ici les arguments. \\

\prop\label{kcor}{Soit $X$ une variété projective et lisse sur $\FF$, de dimension $n$, et soit $l$ un nombre premier, $(l,\mathrm{car}.\FF)=1$. Alors $H^i(X, \hh^{n+1}(\mu_{l^m}^{\otimes j}))=0$ pour $i<n$ et $m>0$.}
\proof{Soit $L/\FF$ une extension galoisienne, de groupe $G$, qui trivialise le faisceau $\mu_{l^m}^{\otimes (j-n)}$ et qui est minimale pour cette propriété. Soit $\pi:Y=X\times_{\FF}L\to X$ le revêtement de $X$ correspondant. Soit $\mathcal F=\hh^{n+1}_Y(\mu_{l^m}^{\otimes j})$ et soit $$0\to \mathcal F\to F^0\to F^1\to\ldots \to F^n\to 0$$ la résolution de Gersten de $\mathcal F$.

\lem\label{acycl}{\begin{itemize} \item [(i)] on a $\hat H^*(G, \pi_*F^i)=0$;
\item [(ii)]$ (\pi_*\mathcal F)^G\simeq \hh^{n+1}_X(\mu_{l^m}^{\otimes j}).$
\end{itemize}}
\proof{ (cf. \cite{Ka93} prop.2) D'après la définition,  $F^i$ est la somme directe sur les points de codimension $i$ de $Y$ des faisceaux $\imath_{y,\,*}H^{n+1-i}(\kappa(y),\mu_{l^m}^{\otimes (j-i)} ))$. Pour établir $(i)$, il suffit donc de voir que $\hat H^*(G,H^{n+1-i}(\kappa(x)\otimes_X{Y},\mu_{l^m}^{\otimes (j-i)} ))=0$ pour $i>0$ et pour tout point $x\in X^{(i)}$. Si $y$ est un point de $Y$ au-dessus de $x$ et $D_y\subset G$ son groupe de décomposition, on a $\hat H^*(G,H^{n+1-i}(\kappa(x)\otimes_X{Y},\mu_{l^m}^{\otimes (j-i)} ))=\hat H^*(D_y,H^{n+1-i}(\kappa(y),\mu_{l^m}^{\otimes (j-i)} ))$ d'après le lemme de Shapiro. On peut donc supposer que $G$ est le groupe de décomposition de $y$.

Sous la condition $\mathrm{cd}_l\,\kappa(x)\leq n+1-i$, on a un isomorphisme
  $$H^{n+1-i}(\kappa(y),\mu_{l^m}^{\otimes (j-i)})_G\to H^{n+1-i}(\kappa(x),\mu_{l^m}^{\otimes (j-i)}),$$ induit par la corestrection. On peut le voir par exemple en utilisant une suite spectrale de Tate (\cite{S}, appendice au ch. 1, Th. 1)
$$E^2_{pq}=H_p(G, H^{n+1-i-q}(\kappa(y), \mu_{l^m}^{\otimes (j-i)} ))\Rightarrow H^{n+1-i-p-q}(\kappa(x), \mu_{l^m}^{\otimes (j-i)}).$$

Puisque $G$ est cyclique ($\FF$ est un corps fini), pour conclure, il suffit d'établir que
 la restriction induit un isomorphisme $H^{n+1-i}(\kappa(x),\mu_{l^m}^{\otimes (j-i)})\to H^{n+1-i}(\kappa(x),\mu_{l^m}^{\otimes (j-i)})^G$. Sous la conjecture de Bloch-Kato, c'est l'énoncé (1) du Th. 1 de Bruno Kahn \cite{Ka93b}.
\qed\\}

Pour montrer la proposition, on utilise la cohomologie mixte de Grothendieck. On a deux suites spectrales (\cite{G57}, 5.2.1):
\begin{gather*}
E_2^{pq}=H^p(G,H^q(X_{Zar},\pi_*\mathcal F))\Rightarrow H^{p+q}(X_{Zar}, G;\pi_*\mathcal F)\\%\label{ss1}
E_2^{'pq}=H^p(X_{Zar},H^q(G,\pi_*\mathcal F))\Rightarrow H^{p+q}(X_{Zar}, G;\pi_*\mathcal F)%\label{s2}
\end{gather*}
où les $H^{i}(X_{Zar}, G; \mathcal G)$ désignent les foncteurs dérivés du foncteur
$$\mathcal G \mapsto \Gamma(X,\mathcal G)^G$$
de la catégorie des faisceaux de Zariski sur $X$, munis d'une action de $G$, compatible avec l'action triviale de $G$ sur $X$, dans la catégorie des groupes abéliens.

D'après le lemme de Shapiro, les termes de la première suite $E_2^{pq}$ se récrivent comme $E_2^{pq}=H^p(G,H^q(Y_{Zar},\mathcal F))$.
Ainsi, pour $p<n$, le théorème \ref{tks} et la suite $E_2^{pq}$ donnent  \begin{equation} H^{p}(X_{Zar}, G;\pi_*\mathcal F)=0.\label{anul}\end{equation}

D'après \cite{Ka93}, prop. 3,  la suite $$0\to H^q(G,\pi_*F^0)\to H^q(G,\pi_*F^1)\to \ldots \to H^q(G,\pi_* F^n)\to 0$$ est une résolution flasque de $H^q(G,\pi_*\mathcal F)$. Ainsi, $H^p(X_{Zar},H^q(G,\pi_*\mathcal F))=0$, $q>0$, d'après le lemme précédent.  La suite $E_2^{'pq}$ et le $(ii)$ du lemme donnent alors des isomorphismes pour tout $p\geq 0$ :
\footnotesize
  \begin{equation*}
H^{p}(X_{Zar}, G;\pi_*\mathcal F)\stackrel{\sim}{\to}H^{p}(X_{Zar}, H^0(G,\pi_*\mathcal F)). \\
  \end{equation*}
\normalsize

D'après \ref{acycl}(ii) et (\ref{anul}), cela donne, pour $0\leq p<n$
\begin{equation*}
 0=H^{p}(X_{Zar}, G;\pi_*\mathcal F)\stackrel{\sim}{\to} H^{p}(X, \hh^{n+1}_X(\mu_{l^m}^{\otimes j})).
\end{equation*}
Cela termine la preuve de la proposition.
 \qed\\}

\subsection{Les invariants}\label{ibf}

%Soit  $\FF$ un corps fini.
%D'après le théorème de Kerz et Saito \cite{KeS}, pour $X$ une $\FF$-variété projective et lisse, géométriquement intègre, et pour $l$ premier différent de $\mathrm{car.}\,\mathbb F$, les groupes $H^i(X, \hh^{n+1}(\QZl(n)))=KH_{n-i}(X, \QZl)$ sont nuls, sauf  $H^n(X, \hh^{n+1}(\QZl(n)))$ qui vaut $\QZl$
Sur un corps fini, on a des analogues du théorème \ref{inv1}.% admet un analogue pour les groupes   $H^i(X, \hh^{n}(\QZl(j)))$.

\theo\label{invf}{Soit  $\FF$ un corps fini. Les groupes $H^i(X, \hh^{n}(\QZl(n-1)))$ où $i<n-1$  et $(l,\mathrm{car.}\FF)=1$, sont des invariants birationnels des $\FF$-variétés projectives et lisses, géométriquement intègres,  de dimension $n$.  De plus,
\begin{itemize}
\item [(i)] pour $X/\FF$ projective et lisse, de dimension $n$ et pour $i<n-1$ on a $H^{i}(X, \hh^{n}(\QZl(n-1)))\simeq H^{i+N}(\mathbb P^N_X, \hh^{n+N}(\QZl(n-1)))$;
\item [(ii)] on a $H^i(\mathbb P^N_k, \hh^{N}(\QZl(n-1)))=0$ pour $i<N-1$;
\item [(iii)] si $X$ et $Y$ sont des $k$-variétés projectives et lisses, géométriquement intègres, de dimensions respectives $n_X\leq n_Y$, stablement birationnelles, i.e. telles que $X\times \mathbb P^{N_X}_k$ et $Y\times \mathbb P^{N_Y}_k$ sont birationnelles, alors $$H^i(X, \hh^{n_X}(\QZl(n-1)))\simeq H^{i+n_Y-n_X}(Y, \hh^{n_Y}(\QZl(n-1))), \quad i< n_X-1.$$
\end{itemize}
\proof{On procède  comme dans la preuve de \ref{inv1}. Avec les mêmes notations, pour établir la première partie, il suffit de montrer que la correspondance $[W]$ agit trivialement sur $H^{i}(X, \hh^{n}(\QZl(n-1)))$. Cela résulte du lemme \ref{lnul}, car les groupes $H^i(V, \hh^{m}(\QZl(n-1)))$, $i<n-1$, sont nuls pour toute $k$-variété projective et lisse $V$, de dimension $n-1$, d'après le théorème de Kerz et Saito \ref{tks}. On obtient l'énoncé $(i)$ en utilisant la suite de localisation et le fait que $H^{i}(X\times\mathbb A^1, \hh^{n+1}(\QZl(n-1)))\stackrel{pr_X^*}{\simeq} H^{i}(X, \hh^{n+1}(\QZl(n-1)))=0$, $i<n-1$, d'après \ref{cks}; $(ii)$ et $(iii)$ en résultent.

\qed\\}

Notons que pour $i=n-1$ le groupe $H^{n-1}(X,\hh^{n}(\QZl(n-1)))$ n'est pas un invariant birationnel. En effet, ce groupe est dual de $H^2_{\acute{e}t}(X, \mathbb Z_l(1))$ (cf.\cite{Mi2} 1.14).\\

Pour des variétés projectives et lisses, de dimension $n$ sur un corps fini $\FF$ on a ainsi le schéma suivant de la page $E_2$ de suite spectrale de Bloch et Ogus \cite{BO} $E_2^{pq}=H^p(X,\hh^{n}(\QZl(n-1)))\Rightarrow H^{p+q}_{\acute{e}t}(X,\QZl(n-1))$, où les lignes $p=0$ et $q=n$ (sauf pour les points $p=n-1$ et $p=n$) consistent d'invariants birationnels (cf. \cite{CT92} 4.4.1 pour la ligne $p=0$).\\

\begin{center}
\includegraphics[scale=0.35]{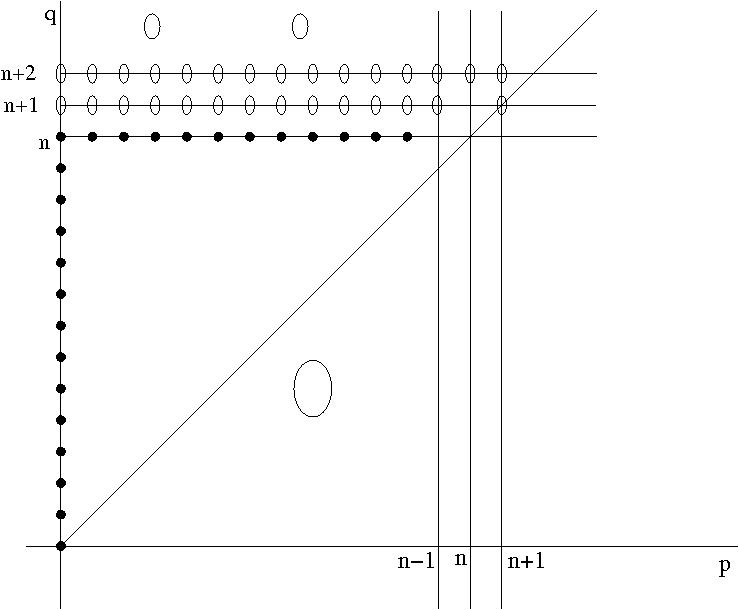}
\end{center}
\vspace{0.5cm}

\rem{L'énoncé du théorème ci-dessus vaut aussi dans les deux cas suivants :
\begin{itemize}
\item[(i)] Pour  les groupes $H^i(X, \hh^{n}(\QZl(j)))$ où $i\geq 0$,  $j\neq n-1$ et $(l,\mathrm{car.}\FF)=1$. Dans ce cas, comme dans la preuve du théorème ci-dessus, pour appliquer le lemme \ref{lnul}, on utilise \ref{cks} pour voir que les groupes $H^i(V, \hh^{n}(\QZl(n-1)))$ sont nuls pour $V$ projective et lisse de dimension $n-1$. Cela n'utilise donc pas la conjecture de Kato.
\item[(ii)] Pour les groupes  $H^i(X, \hh^{n+1}(\mu_{l^m}^{\otimes j}))$ où $i<n-1$ et $(l,\mathrm{car.}\FF)=1$. Dans ce cas on utilise \ref{kcor} pour appliquer  le lemme \ref{lnul}.\\
\end{itemize}}

\subsection{Lien avec les $1$-cycles}\label{luc}

Pour $k$ un corps, $X$ un $k$-schéma lisse et $A$ un groupe abélien, on dispose d'un complexe de faisceaux $A_X(n)_{\acute{e}t}$ (resp. $A_X(n)$ pour le complexe de faisceaux Zariski), défini à partir du complexe  de cycles de Bloch $z^n(X,\cdot)$ (cf. par exemple  \cite{Ka11} 2.2). On note $\HH^{i}(X,A_X(n))$  (resp. $\mathbb H^{i}(X,A_X(n))$ ) l'hypercohomologie de ce complexe. Pour $n=0$ on a $\mathbb Z_X(0)=\mathbb Z$. Pour $n=1$ le complexe  $\mathbb Z_X(1)$ est quasi-isomorphe à $\mathbb G_m[-1]$ (cf. par exemple  \cite{Ka11} (2.4)).

Pour $X$ quasi-projectif et lisse, on a $CH^n(X,2n-i)\simeq \mathbb H^{i}(X,\ZZ_X(n))$.

Pour $X$ lisse, on a une suite  spectrale de coniveau ( \cite{Ka10} 2.7), \begin{equation*}
                                                                      E_1^{p,q}=\bigoplus_{x\in X^{(p)}}\HH^{q-p}(\kappa(x),\mathbb Z(n-p))\Rightarrow \HH^{p+q}(X, \mathbb Z_X(n))
                                                                       \end{equation*}
  où
\begin{enumerate}\item pour $q\geq n+2$,  $E_2^{p,q}=H^p(X,\mathcal H^{q-1}(\mathbb Q/\mathbb Z(n)))$, cf. \cite{Ka97} 5.1.
\item $E_1^{p,q}=0$ pour $q=n+1$ (sous la conjecture de Bloch-Kato, qui est maintenant établie par Voevodsky \cite{Voe}), cf. \cite{Ka10} 2.7.\\
\end{enumerate}

\theo\label{im}{\textit{Soit $X$ une $\FF$-variété projective et lisse, de dimension $n$. On a une suite exacte à la $\mathrm{car}.\FF$-torsion près
\footnotesize
 $$H^{n-4}(X,\mathcal H^n(\QZ(n-1)))\to CH^{n-1}(X)\to\HH^{2n-2}(X,\mathbb Z_X(n-1))\to H^{n-3}(X,\mathcal H^n(\QZ(n-1)))\to 0$$
}}
\normalsize
\vspace{-0.6cm}
\proof{On dispose d'une suite spectrale de coniveau pour  $\mathbb Z_X(n-1))$ : $$E_1^{p,q}=\bigoplus_{x\in X^{(p)}}\HH^{q-p}(\kappa(x),\mathbb Z(n-1-p))\Rightarrow \HH^{p+q}(X, \mathbb Z_X(n-1)).  $$ On a :
\begin{enumerate}
 \item  $E_2^{p,q}= E_1^{p,q}=0$ pour $p\geq n+1$ car $ X^{(p)}$ est vide pour un tel $p$.
\item  la ligne $E_2^{p,n}$ est constituée de zéros, sous la conjecture de Bloch-Kato;
\item  $E_2^{p,n+1}=H^p(X,\mathcal H^{n}(\mathbb Q/\mathbb Z(n-1)))$;
\item à la $p$-torsion près, $E_2^{p,n+2}=H^p(X,\mathcal H^{n+1}(\mathbb Q/\mathbb Z(n-1)))=0$  d'après \ref{cks};
\item pour $q\geq n+3$, $E_2^{p,q}=H^p(X,\mathcal H^{q-1}(\mathbb Q/\mathbb Z(n)))=0$ d'après Bloch-Ogus ($dim\,X=n$ et $cd\,\FF\leq 1$);
\item pour $q=n-1$ et pour $q=n-2$, $E_2^{n,q}= E_1^{n,q}=0$. En effet, pour $F$ un corps et pour $m<0$, on a   $H^r(F,\mathbb Z(m))= H^{r-1}(F,\QZ(m-1))=0$ pour $r\leq 0$, cf. \cite{Ka97} 2.4.
\item $E_2^{n-1,n-1}=CH^{n-1}(X)$. En effet,   $E_1^{n-1,n-1}=\bigoplus_{x\in X^{(n-1)}} \HH^{0}(\kappa(x),\mathbb Z)=\bigoplus_{x\in X^{(n-1)}} \mathbb Z$ et $E_1^{n-2,n-1}=\bigoplus_{x\in X^{(n-2)}} \HH^{1}(\kappa(x),\mathbb Z(1))=\bigoplus_{x\in X^{(n-2)}} \kappa(x)^*$.\\% La flèche $E_1^{n-2,n-1}\to E_1^{n-1,n-1}$ provient de est donnée par des flèches diviseurs
\end{enumerate}

Ainsi, à la $p$-torsion près, dans la filtration sur $\HH^{2n-2}(X,\mathbb Z_X(n-1))$, on n'a que deux termes  $$E_{\infty}^{n-3,n+1}=E_{2}^{n-3,n+1}=H^{n-3}(X,\mathcal H^{n}(\mathbb Q/\mathbb Z(n-1)))$$ et  $$E_{\infty}^{n-1,n-1}=E_2^{n-1,n-1}/E_2^{n-4,n+1}=CH^{n-1}(X)/H^{n-4}(X,\mathcal H^n(\QZ(n-1))).$$ On obtient ainsi la suite du théorème. \qed\\
}

La méthode de Colliot-Thélène et Kahn  \cite{CTK} 2.1 permet de relier le groupe $H^{n-3}(X,\mathcal H^n(\QZl(n-1)))$ avec le conoyau  de l'application classe de cycle
$CH^{n-1}(X)\otimes \ZZl \to H^{2n-2}(X,\ZZl(n-1))$. Rappelons que pour un corps $k$ on dit que $k$ est à cohomologie galoisienne finie, si pour tout module galoisien fini $M$ et tout entier $i\geq 0$, les groupes de cohomologie $H^i(k,M)$ sont finis. L'énoncé suivant pour $i=2$ est dans \cite{CTK} 2.1. La preuve du cas général utilise les mêmes arguments, on les détaille ci-dessous.

\theosd\label{me}(cf. \cite{CTK} 2.1) {\textit{ Soient $k$ un corps à cohomologie galoisienne finie, $X$ une $k$-variété lisse et $l$ un nombre premier, $(l,\mathrm{car}.k)=1$. Soit $i\geq 0$. Posons
 $$M=Coker [CH^{i}(X)\otimes \ZZl \to H^{2i}_{\acute{e}t}(X,\ZZl(i))]$$ le conoyau de l'application classe de cycle $l$-adique étale
et
$$C=Coker [CH^{i}(X) \to \HH^{2i}(X,\mathbb Z_X(i))]$$  le conoyau de l'application classe de cycle motivique étale.\\
Alors le groupe (fini) de torsion $M_{tors}$ est  isomorphe au quotient du groupe $C\{l\}$ de torsion $l$-primaire de $C$ par son sous-groupe divisible maximal.\\
}

\rem\label{rt}{Le groupe $C$ est de torsion. En effet, il résulte de \cite{Ka11} (2.2) et 2.6(c) qu'on a un isomorphisme $CH^{i}(X)\otimes  \mathbb Q\stackrel{\sim}{\to} H^{2i}_{Zar} (X, \mathbb Q_X(i))\stackrel{\sim}{\to}\HH^{2i}(X, \mathbb Q_X(i))$. \\}

En combinant le théorème ci-dessus avec le théorème \ref{im}, on obtient

\corro{\textit{Soit $X$ une $\FF$-variété projective et lisse, de dimension $n$ et soit $l$ un nombre premier, $(l, \mathrm{car}.\FF)=1$.  Alors le groupe (fini) de torsion du  conoyau  de l'application classe de cycle
$CH^{n-1}(X)\otimes \ZZl \to H^{2n-2}_{\acute{e}t}(X,\ZZl(n-1))$ est isomorphe au quotient du groupe $H^{n-3}(X,\mathcal H^n(\QZl(n-1)))$ par son sous-groupe divisible maximal. \\}}

\rem{\begin{itemize}
\item[(i)] Si $n\geq 2$ et si  $\mathrm{Br}\,X$ est fini (ce qui est le cas sous la conjecture de Tate pour les cycles de codimension $1$), le conoyau de l'application $CH^{n-1}(X)\otimes \ZZl \to H^{2n-2}(X,\ZZl(n-1))$ est fini %et nul pour presque tout $l$
(\cite{CTK} 7.3).
      \item [(ii)] Pour $X$ une variété appartenant à la classe $B_{\mathrm{Tate}}(\FF)$ (conjecturalement, c'est le cas pour toute  $\FF$-variété projective et lisse), le groupe $\HH^{2n-2}(X,\mathbb Z_X(n-1))$ est de type fini (\cite{Ka03} 3.10). 
Ainsi, le groupe $H^{n-3}(X,\mathcal H^n(\QZl(n-1)))$ est conjecturalement fini.   \\
     \end{itemize}
}

\rem{ Soit $X$ une $\FF$-variété projective et lisse, de dimension $n$, munie d'un morphisme $\phi :X\to C$ où $C$ est une courbe projective et lisse, géométriquement connexe. Soit $l$ un nombre premier, $(l, \mathrm{car}.\FF)=1$. D'après un théorème de Saito \cite{Sai}, si l'application classe de cycle $l$-adique
$CH^{n-1}(X)\otimes \ZZl\to H^{2n-2}_{\acute{e}t}(X,\ZZl(n-1))$ est surjective, alors  l'existence d'un zéro-cycle de degré premier à $l$ sur  la fibre générique $X_{\eta}$ de l'application $\phi$ se détecte seulement en utilisant  l'obstruction de Brauer-Manin. Plus précisement, s'il existe sur  $X_{\eta}$ une famille de zéro-cycles locaux de degré $1$, orthogonale au groupe de Brauer de $X_{\eta}$ via l'accouplement de Brauer-Manin (cf. \cite {Man71}), alors il existe sur $X_{\eta}$ un zéro-cycle de degré premier à $l$. }
$$\,$$

\textit{Démonstration du théorème \ref{me}.}
Pour montrer ce théorème, on va d'abord établir qu'on a une application naturelle $\phi: \HH^{2i}(X,\mathbb Z_X(i))\otimes\ZZl\to H^{2i}(X,\ZZl(i))$ telle que
\begin{itemize}
 \item [(i)]  conoyau de $\phi$ est sans torsion;
\item[(ii)] l'application induite $$\phi/l^r: [\HH^{2i}(X,\mathbb Z_X(i))\otimes\ZZl]/l^r\to H^{2i}(X,\ZZl(i))/l^r$$ est injective.
\end{itemize}

On dispose d'un triangle exact (cf. \cite{Ka11} 2.6)
$$\mathbb Z_{\acute{e}t}(i)\stackrel{\times l^r}{\to}\mathbb Z_{\acute{e}t}(i)\to \mu_{l^r}^{\otimes i}.$$
En prenant l'hypercohomologie, on en déduit une suite exacte
\begin{equation}\label{se}
 0\to \HH^{2i}(X,\mathbb Z_X(i))/l^r \to H^{2i}_{\acute{e}t}(X,\mu_{l^r}^{\otimes i})\to \HH^{2i+1}(X,\mathbb Z_X(i))[l^r]\to 0.
\end{equation}
Puisque le groupe de milieu est fini, on obtient une suite exacte en passant à la limite projective:
\begin{equation}\label{st} 0\to \varprojlim [\HH^{2i}(X,\mathbb Z_X(i))/l^r] \to H^{2i}_{\acute{e}t}(X,\mathbb Z_l(i))\to \varprojlim[\HH^{2i+1}(X,\mathbb Z_X(i))[l^r]]\to 0.
\end{equation}

Les applications $\HH^{2i}(X,\mathbb Z_X(i))\to \HH^{2i}(X,\mathbb Z_X(i))/l^r $ induisent une application  $\HH^{2i}(X,\mathbb Z_X(i))\otimes\ZZl\to \varprojlim [\HH^{2i}(X,\mathbb Z_X(i))/l^r]$ qui est surjective d'après \cite{CTK} 2.2 (ii). On définit alors $\phi$ comme  une application composée $$\phi:\HH^{2i}(X,\mathbb Z_X(i))\otimes\ZZl\to \varprojlim [\HH^{2i}(X,\mathbb Z_X(i))/l^r]\to H^{2i}_{\acute{e}t}(X,\mathbb Z_l(i)).$$ Son conoyau est sans torsion car le terme de droite de la suite (\ref{st}) l'est, puisque c'est un module de Tate.

L'injectivité de $\phi/l^r$ résulte du diagramme commutatif suivant
\footnotesize
$$\begin{CD}  [\HH^{2i}(X,\mathbb Z_X(i))\otimes\ZZl]/l^n@>>>    \HH^{2i}(X,\mathbb Z(i))/l^n.\\
@V\phi/l^n VV  @VVV @. \\
 H^{2i}_{\acute{e}t}(X,\mathbb Z_l(i))/l^n @>>>   H^{2i}_{\acute{e}t}(X,\mu_{l^r}^{\otimes i}).
\end{CD}$$
\normalsize
où les applications du bas et de droite sont injectives (cf.  la suite (\ref{se}) pour celle de droite) et l'application  du haut est un isomorphisme d'après \cite{CTK} 2.2 (i).

Puisque $C$ est de torsion (cf. \ref{rt}), on a un diagramme commutatif suivant

\footnotesize
$$\begin{CD}  CH^{i}(X)\otimes \ZZl@>>> \HH^{2i}(X,\mathbb Z_X(i))\otimes\ZZl@>>>  C\{l\}@>>> 0\\
@V\simeq VV @V\phi VV  @V\psi_0 VV @. \\
CH^{i}(X)\otimes \ZZl@>>>  H^{2i}(X,\ZZl(i))@>>>  M@>>> 0.
\end{CD}$$
\normalsize
où l'application $\psi_0$ est induite par le diagramme. Notons que
la commutativité du carré de gauche de ce diagramme, i.e. la compatibilité des applications classe de cycle résulte de \cite{Ka11} 3.1.

Comme  $C\{l\}$ est un groupe de torsion, l'application $\psi_0$ induit une application $\psi: C\{l\}\to M_{tors}$. Puisque le conoyau de $\phi$ est sans torsion, on déduit par chasse au diagramme que $\psi$ est surjective.
Puisqu'on a de plus que $\phi/l^n$ est injective, on en déduit qu'on a un isomorphisme $\psi/l^n: C\{l\}/l^n\to M_{tors}/l^n$.

Puisque le groupe $M_{tors}$ est fini, on a $M\simeq M/l^N$ pour $N$ assez grand, i.e. pour $N\geq N_0$.  Pour un tel $N$, on a donc un isomorphisme $l^{N+1}C\{l\}\simeq l^NC\{l\}$. Ainsi, le groupe $l^{N_0}C\{l\}$ est le sous-groupe divisible maximal de $C\{l\}$ et le quotient de $C\{l\}$ par ce sous-groupe et isomorphe à $M_{tors}$.\qed

\end{document}